\documentclass[12pt,english,a4paper]{amsart}
\usepackage[pass]{geometry}
\usepackage[cp1251]{inputenc}
\usepackage[english]{babel}
\usepackage{amsmath,amsthm,amssymb}
\newtheorem{theorem}{Theorem}[section]
\newtheorem{lemma}[theorem]{Lemma}
\newtheorem{corollary}[theorem]{Corollary}

\numberwithin{equation}{section}

\theoremstyle{definition}

\newtheorem{example}[theorem]{Example}

\title{Empirical approximation of $L_p$ norms for $p\in (2, \infty)$}

\author{Feng Dai}
\address{F.~Dai, Department of Mathematical and Statistical Sciences, University of Alberta, Edmonton, Alberta T6G 2G1, Canada}
\email{fdai@ualberta.ca}

\author{Egor Kosov}
\address{E.~Kosov, Centre de Recerca Matem\`atica, Campus de Bellaterra, Edifici~C 08193
Bellaterra (Barcelona), Spain.}
\email{kosoved09@gmail.com}

\author{Noel Murasko}
\address{N.~Murasko, Department of Mathematical and Statistical Sciences, University of Alberta, Edmonton, Alberta T6G 2G1, Canada}
\email{murasko@ualberta.ca}

\begin{document}

\begin{abstract}
For $p>2$ and an $N$-dimensional subspace
$X\subset L_\infty(\mu)$
satisfying
\[
\|f\|_\infty\le (KN)^{1/p}\|f\|_{L_p(\mu)}
\quad\forall f\in X,
\]
we show that
$m\ge C(p)KN
\bigl(\varepsilon^{-p}
+\varepsilon^{-2}\log_2(2N)\bigr)$
independent samples $\xi^1,\ldots,\xi^m$ with distribution $\mu$
satisfy
\[
(1-\varepsilon)\|f\|_{L_p(\mu)}^p
\le\frac1m\sum_{j=1}^m|f(\xi^j)|^p
\le(1+\varepsilon)\|f\|_{L_p(\mu)}^p
\quad\forall f\in X
\]
with probability at least $1-\exp(-c(p)\frac{m\varepsilon^2}{KN})$.
This gives the optimal sample-size order
$N\log_2 N$ for independent random sampling with fixed
$K$, $p$, and $\varepsilon$.
This result follows from a more general empirical moment
estimate over unit balls of $p$-convex norms.
\end{abstract}
\maketitle

\section{Introduction}

For $p\ge2$, a norm $\|\cdot\|$ on a linear space $X$ is called $p$-convex
with constant $R>0$ if
\[
\biggl\|\frac{x+y}{2}\biggr\|
\le
\max\{\|x\|,\|y\|\}
-R\|x-y\|^p
\quad\forall x,y\in X\text{ with }\|x\|,\|y\|\le1.
\]

In this paper we prove the following theorem.

\begin{theorem}\label{th-sampling-p-convex}
Let $p>2$, $H\ge M>0$, let $(\Omega,\mu)$ be a probability space, and let
$X\subset L_\infty(\mu)$ be an $N$-dimensional real
Banach space, whose norm $\|\cdot\|$ is $p$-convex with constant $R$.
Assume that
\[
\sup\limits_{\|f\|\le 1}\|f\|_{\infty}^p\le H
\quad\text{and}\quad \sup\limits_{\|f\|\le 1}\|f\|_{L_p(\mu)}^p\le M.
\]
Then there exist constants $C(p,R)>0$ and $c(p)>0$ such that,
for every $\varepsilon\in(0,1)$ and every integer
\begin{equation}\label{eq-sample-complexity}
m\ge C(p,R)\frac{H}{M}
\Bigl(\varepsilon^{-p}
+\varepsilon^{-2}\log_2(2N)\Bigr),
\end{equation}
independent samples $\xi^1,\ldots,\xi^m$ with distribution $\mu$
satisfy
\[
\biggl|
\frac1m\sum_{j=1}^m|f(\xi^j)|^p
-\|f\|_{L_p(\mu)}^p
\biggr|
\le \varepsilon M\|f\|^p
\quad\forall f\in X
\]
with probability at least
$1-\exp\bigl(-c(p)\frac{m\varepsilon^2M}{H}\bigr)$.
\end{theorem}

\subsection{Consequences and applications}

Theorem~\ref{th-sampling-p-convex} can be equivalently formulated
in terms of empirical approximation of the $p$-th moments
of one-dimensional marginals of random vectors
uniformly over unit balls of $p$-convex norms.

\begin{corollary}\label{cor-random-vectors-p-convex}
Let $p>2$, $H\ge M>0$, let ${\bf Z}$ be a random vector in $\mathbb R^N$,
and let $V\subset\mathbb R^N$ be the unit ball of a
$p$-convex norm with constant $R$.
Assume that
\[
\|{\bf Z}\|_{V^\circ}\le H^{1/p}
\quad\text{almost surely}
\quad\text{and}\quad
\sup_{y\in V}\mathbb E|\langle{\bf Z},y\rangle|^p\le
M.
\]

There exist constants $C(p,R)>0$ and $c(p)>0$ such that,
for every $\varepsilon\in(0,1)$ and every integer
\[
m\ge C(p,R)\frac HM
\bigl(\varepsilon^{-p}
+\varepsilon^{-2}\log_2(2N)\bigr),
\]
independent copies ${\bf Z}_1,\ldots,{\bf Z}_m$ of ${\bf Z}$
satisfy
\[
\sup_{y\in V}
\biggl|
\frac1m\sum_{j=1}^m|\langle{\bf Z}_j,y\rangle|^p
-\mathbb E|\langle{\bf Z},y\rangle|^p
\biggr|
\le\varepsilon M
\]
with probability at least
$1-\exp\bigl(-c(p)\frac{m\varepsilon^2M}{H}\bigr)$.
\end{corollary}

Gu\'edon and Rudelson~\cite{GR} studied such empirical
approximation of moments over $p$-convex sets under the
additional assumption $V\subset DB_2^N$.
Their result gives the same uniform approximation provided that
$\|{\bf Z}\|_2\le L$ almost surely and
\[
m\ge C(p,R)\varepsilon^{-2}\frac{(DL)^p}{M}
\Bigl[\log_2\Bigl(\frac{2(DL)^p}{M\varepsilon}\Bigr)\Bigr]^{2-2/p}.
\]
This sufficient condition was improved in
\cite[Corollary~4.4]{Kosov2021} to
\[
m\ge C(p,R)\varepsilon^{-2}
\frac{(DL)^p}{M}
\log_2\Bigl(\frac{2(DL)^p}{M\varepsilon}\Bigr).
\]
Our more general result also applies in this special setting,
where it recovers this bound with an improved logarithmic factor
and a worse dependence on $\varepsilon$.
Indeed, for $V\subset DB_2^N$, one may take $H=(DL)^p$ in
Corollary~\ref{cor-random-vectors-p-convex}, which gives
\[
m\ge C(p,R)\frac{(DL)^p}{M}
\Bigl(\varepsilon^{-p}
+\varepsilon^{-2}\log_2(2N)\Bigr).
\]

Applying Theorem~\ref{th-sampling-p-convex} to finite-dimensional
subspaces equipped with the $L_p(\mu)$ norm, we obtain the
following sampling discretization result for $L_p$ norms when $p>2$.

\begin{corollary}\label{cor-sampling-Lp}
Let $p>2$, $K\ge1$, let $(\Omega,\mu)$ be a probability space, and let
$X\subset L_\infty(\mu)$ be a finite-dimensional real
or complex linear subspace of dimension $N$. Suppose that
\[
\|f\|_{\infty}
\le (KN)^{1/p}\|f\|_{L_p(\mu)}
\quad\forall f\in X.
\]
Then there exist constants $C(p)>0$ and $c(p)>0$ such that,
for every $\varepsilon\in(0,1)$ and every integer
\[
m\ge C(p)KN
\bigl(\varepsilon^{-p}
+\varepsilon^{-2}\log_2(2N)\bigr),
\]
independent samples $\xi^1,\ldots,\xi^m$ with distribution $\mu$
satisfy
\[
(1-\varepsilon)\|f\|_{L_p(\mu)}^p
\le
\frac1m\sum_{j=1}^m|f(\xi^j)|^p
\le
(1+\varepsilon)\|f\|_{L_p(\mu)}^p
\quad\forall f\in X
\]
with probability at least
$1-\exp\bigl(-c(p)\frac{m\varepsilon^2}{KN}\bigr)$.
\end{corollary}

This result is essentially optimal for independent random sampling,
as illustrated by the examples in Section~\ref{sect-examples}.
It improves the previous sufficient bound of
$C(p,\varepsilon)KN(\log_2(KN))^{\min\{3,p\}}$
sampling points under the same assumption, obtained in
\cite[Corollary~4.8]{Kosov2021} and
\cite[Theorem~2.2]{DT1}.
See the surveys \cite{DKT26,KKLT}
for further discussion of this and related problems.

\subsection{Strategy of the proof and organization of the paper}

The main step in the proof of Theorem~\ref{th-sampling-p-convex}
is an estimate for sums of entropy numbers of linear operators
defined on Banach spaces whose duals have type $q$.

For a bounded linear operator $L\colon X\to Y$ and $j\in\mathbb N$,
let
\[
\varepsilon_j(L)
:=
\inf\Bigl\{
\rho>0\colon
LB_X\subset\bigcup_{k=1}^{2^{j-1}}({\bf y}_k+\rho B_Y),\
{\bf y}_k\in Y
\Bigr\},
\]
where $B_X$ and $B_Y$ denote the unit balls of $X$ and $Y$,
respectively. Let
\[
\mathcal E_{2,q}(L)
:=
\biggl(
\sum_{j=1}^\infty
j^{q/2-1}\bigl[\varepsilon_j(L)\bigr]^q
\biggr)^{1/q}.
\]

Recall that, for $1\le q\le2$, a Banach space $Y$ has
Rademacher type $q$ with constant $T$ if, for every finite
family ${\bf y}_1,\ldots,{\bf y}_n\in Y$,
\[
\mathbb E
\Bigl\|\sum_{j=1}^n\eta_j{\bf y}_j\Bigr\|_Y^q
\le T^q\sum_{j=1}^n\|{\bf y}_j\|_Y^q,
\]
where $\eta_1,\ldots,\eta_n$ are independent Rademacher
random variables.

We prove the following estimate.

\begin{theorem}\label{th-row-enrtopy}
Let $1<q<2$, and let $X$ be a
Banach space whose dual has type $q$ with constant at most $T$.
For $a_1,\ldots,a_m\in X^*$, consider
$A\colon X\to\ell_2^m$,
$Ax:=(a_1(x),\ldots, a_m(x))$.
Then
\[
\mathcal E_{2,q}(A)
\le C(q,T)
\biggl(\sum_{k=1}^m\|a_k\|_{X^*}^q\biggr)^{1/q}.
\]
\end{theorem}

By entropy duality,
this estimate follows from its dual formulation,
recorded in Theorem~\ref{th-column-entropy},
which is a quantitative consequence of
a result of Cobos \cite[Theorem~4.1]{Cobos}.

We postpone the proof of Theorem~\ref{th-row-enrtopy} to
Section~\ref{sect-entropy}. First, in Section~\ref{sect-discr},
we deduce Theorem~\ref{th-sampling-p-convex}
and Corollary \ref{cor-sampling-Lp} from
Theorem~\ref{th-row-enrtopy}.
Symmetrization and Talagrand's chaining theorem reduce the sampling
problem to bounds for the chaining functional
$\gamma_{2,1}(B_E,\varrho_p)$, where $B_E$ is the unit ball of a
$p$-convex norm on a subspace $E\subset\mathbb R^m$ and
\[
\varrho_p(x,y)
:=
\biggl(
\sum_{j=1}^m
|x_j-y_j|^2
\bigl(|x_j|^{p-1}+|y_j|^{p-1}\bigr)^2
\biggr)^{1/2}.
\]
We apply the chaining reduction from
\cite[Theorem~2.6 and Lemma~3.3]{Kosov2021},
which estimates $\gamma_{2,1}(B_E,\varrho_p)$ by the sum of
\[
\bigl[\gamma_{2p,p}(B_E,\|\cdot\|_{2p})\bigr]^p
\quad
\text{and}
\quad
\sup_{z\in B_E}
\biggl(
\sum_{n=0}^\infty
\bigl[
2^{n/2}e_n(B_E,\|\cdot\|_z)
\bigr]^q
\biggr)^{1/q},
\]
where $q=\frac{p}{p-1}$ and
$\|x\|_z
:=
\bigl(
\sum_{k=1}^m|z_k|^{2p-2}|x_k|^2
\bigr)^{1/2}$.
For the first term, we use the known estimate from
\cite[Lemma~3.4(3) and Corollary~4.2]{Kosov2021}.
To control the second, we apply Theorem~\ref{th-row-enrtopy}.

Finally, Section~\ref{sect-examples} contains examples showing
that the sampling discretization result in
Corollary~\ref{cor-sampling-Lp} is essentially optimal.

\subsection{Notation}
For $1\le p<\infty$, we equip $\mathbb R^m$ with the norm
\[
\|x\|_p:=\biggl(\sum_{j=1}^m|x_j|^p\biggr)^{1/p}.
\]
For a normed space $X$, let
\[
B_X:=\{x\in X\colon \|x\|_X\le1\}.
\]
Let
$\mathbb N_0:=\mathbb N\cup\{0\}$.
For a nonempty bounded subset $A$ of a normed space $X$ and
$n\in\mathbb N_0$, let
\[
e_n(A,\|\cdot\|_X)
:=
\inf_{\substack{\Lambda\subset A\\1\le|\Lambda|<2^{2^n}}}
\sup_{x\in A}\min_{y\in\Lambda}\|x-y\|_X.
\]

Let $(T,\varrho)$ be a nonempty metric or quasi-metric space.
A sequence $(\mathcal A_n)_{n\ge0}$ of partitions of $T$
is called admissible if $\mathcal A_{n+1}$ refines
$\mathcal A_n$,
\[
\mathcal A_0=\{T\},
\quad
|\mathcal A_n|<2^{2^n}\quad\forall n\ge1.
\]
For $t\in T$, let $\mathcal A_n(t)$ denote the unique
element of $\mathcal A_n$ containing $t$.
For $\alpha>0$ and $\beta\ge1$, define the chaining functional
\[
\gamma_{\alpha,\beta}(T,\varrho)
:=
\inf_{(\mathcal A_n)}
\sup_{t\in T}
\biggl(
\sum_{n=0}^\infty
\bigl[
2^{n/\alpha}
\operatorname{diam}(\mathcal A_n(t),\varrho)
\bigr]^\beta
\biggr)^{1/\beta},
\]
where the infimum is taken over all admissible sequences,
and
\[
\operatorname{diam}(A,\varrho):=\sup_{x,y\in A}\varrho(x,y).
\]
When $\varrho(x,y)=\|x-y\|$, we use the notation
$\gamma_{\alpha,\beta}(T,\|\cdot\|)$.

The letters $c,C$ denote positive absolute constants whose
values may change from one occurrence to another.
Constants depending on $p$ are denoted by $c(p)$ and $C(p)$,
and those depending on $p$ and $R$ by $C(p,R)$. When several constants appear in the same argument, we distinguish them by subscripts,
writing $C_1$, $C_2$, and so on. These subscripts are only labels and do not indicate any additional
dependence.

\section{Sampling discretization}
\label{sect-discr}

\subsection{Moment estimate for the power process}

\begin{lemma}\label{lem-signed-sum}
Let $p> 2$, let $E\subset\mathbb R^m$ be a linear subspace,
endowed with a $p$-convex norm $\|\cdot\|$ with constant
$R$,
and assume
\[
\sup_{x\in B_E}\|x\|_\infty^p\le h\quad \text{and}\quad \sup_{x\in B_E}\|x\|_p^p\le M,
\]
where $B_E:= \{x\in E\colon \|x\|\le 1\}$.
Then,
for i.i.d.\! Rademacher random variables $\eta_1,\ldots,\eta_m$ and every $r\ge1$, one has
\[
\Bigl[
\mathbb E
\sup_{x\in B_E}
\Bigl|\sum_{k=1}^m\eta_k |x_k|^p\Bigr|^r
\Bigr]^{1/r}
\le
C(p,R)
\big[h^{1/p}M^{1-1/p}+\sqrt{hM\log_2(2m)}\big]+C(p)\sqrt {rhM}.
\]
\end{lemma}

\begin{proof}
Let
\[
\varrho_p(x,y):=
\Big(
\sum_{j=1}^m
|x_j-y_j|^2
\bigl(|x_j|^{p-1}+|y_j|^{p-1}\bigr)^2
\Big)^{1/2}.
\]
By integrating the tail estimate in \cite[Theorem 2.2.27]{Talagrand2014} and since
\[
\Bigl(\sum_{k=1}^m
\bigl||x_k|^p-|y_k|^p\bigr|^2\Bigr)^{1/2}
\le p\varrho_p(x,y),
\]
we have
\begin{equation}\label{eq-moment-est}
\Bigl[
\mathbb E
\sup_{x\in B_E}
\Bigl|\sum_{k=1}^m\eta_k |x_k|^p\Bigr|^r
\Bigr]^{1/r}
\le
C_1(p)\gamma_{2, 1}(B_E, \varrho_p)+C_2(p)\sqrt{r}\operatorname{diam}(B_E, \varrho_p).
\end{equation}
Since $\varrho_p$ is a quasi-metric,
\[
\varrho_p(x, y)\le C_3(p)\bigl(\varrho_p(x, 0) +
\varrho_p(y, 0)\bigr)
\]
and, for $x\in B_E$, we have
\[
\varrho_p(x, 0) = \|x\|_{2p}^p
\le \|x\|_\infty^{p/2}\|x\|_p^{p/2}
\le \sqrt{hM}.
\]
Therefore,
\begin{equation}\label{eq-diam-est}
\operatorname{diam}(B_E, \varrho_p)
\le 2C_3(p)\sqrt{hM},
\end{equation}
and
it remains to control
$\gamma_{2, 1}(B_E, \varrho_p)$.

By \cite[Theorem~2.6]{Kosov2021} and by the argument in the proof of \cite[Lemma~3.3]{Kosov2021} for $r=2$, we have
\begin{equation}\label{eq-gamma-est}
\gamma_{2,1}(B_E,\varrho_p)
\le C_4(p, R)\bigg[
\sup_{z\in B_E}
\bigg(
\sum_{n=0}^\infty\bigl[2^{n/2}
e_n(B_E,\|\cdot\|_z)\bigr]^q
\bigg)^{1/q}
+\bigl[
\gamma_{2p,p}(B_E,\|\cdot\|_{2p})
\bigr]^p
\bigg],
\end{equation}
where $q=\frac{p}{p-1}$, and,
for $z\in B_E$,
\[
\|x\|_z:=
\bigg(
\sum_{k=1}^m |z_k|^{2p-2}|x_k|^2
\bigg)^{1/2}.
\]
Applying \cite[Lemma 3.4(3)]{Kosov2021} with
$r=2$, $\theta=p$, and $s=p$, we obtain
\[
\bigl[\gamma_{2p,p}(B_E,\|\cdot\|_{2p})\bigr]^p
\le C_5(p,R)\sqrt{M}
\sup_{n\ge0}2^{n/2}\bigl[e_n(B_E,\|\cdot\|_\infty)\bigr]^{p/2}.
\]
By \cite[Corollary 4.2]{Kosov2021},
\[
e_n(B_E,\|\cdot\|_\infty)
\le C_6(p, R)h^{1/p}\,2^{-n/p}\bigl[\log_2(2m)\bigr]^{1/p}.
\]
Consequently,
\begin{equation}\label{eq-2-term}
\bigl[\gamma_{2p,p}(B_E,\|\cdot\|_{2p})\bigr]^p
\le C_7(p,R)\sqrt{hM\log_2(2m)}.
\end{equation}

Let $A\colon E\to \ell_2^m$,
$Ax := (a_1(x), \ldots, a_m(x))$,
where $a_k(x):= |z_k|^{p-1}x_k$.
Since
\[
e_n(B_E,\|\cdot\|_z) \le
2\varepsilon_{2^n}(A)
\]
and due to monotonicity of the entropy numbers,
we have
\[
\bigg(
\sum_{n=0}^\infty\bigl[2^{n/2}
e_n(B_E,\|\cdot\|_z)\bigr]^q
\bigg)^{1/q}
\le
C_8(q)\biggl(\sum_{j=1}^\infty j^{q/2-1}\bigl[\varepsilon_j(A)\bigr]^q\biggr)^{1/q}.
\]
For $p>2$, the dual space $E^*$ has Rademacher type
$q=p/(p-1)$ with constant $T=T(p,R)$, see, for example, \cite[Proposition~1.e.2 and Theorem~1.e.16]{LT1979}.
Theorem \ref{th-row-enrtopy} gives
\[
\biggl(\sum_{j=1}^\infty
j^{q/2-1}[\varepsilon_j(A)]^q\biggr)^{1/q}
\le C_9(q,T)
\biggl(\sum_{k=1}^m\|a_k\|_{E^*}^q\biggr)^{1/q}.
\]
Since
\[
\|a_k\|_{E^*} = |z_k|^{p-1}\sup_{x\in B_E}|x_k|
\le h^{1/p}|z_k|^{p-1},
\]
we obtain
\[
\sup_{z\in B_E}\bigg(
\sum_{n\ge0}\bigl[2^{n/2}
e_n(B_E,\|\cdot\|_z)\bigr]^q
\bigg)^{1/q}
\le C_{10}(p, R)h^{1/p}
M^{1-1/p}.
\]
Combining this estimate with \eqref{eq-2-term} and \eqref{eq-gamma-est}, we conclude
\[
\gamma_{2,1}(B_E,\varrho_p)\le C_{11}(p, R)\bigl(h^{1/p}
M^{1-1/p} + \sqrt{hM \log_2(2m)}\bigr),
\]
which together with \eqref{eq-diam-est} and \eqref{eq-moment-est} implies the assertion.
\end{proof}

\subsection{Proof of Theorem \ref{th-sampling-p-convex}}
We first prove the following lemma, which allows us to obtain the sample-size bound \eqref{eq-sample-complexity} with a logarithmic factor depending only on $N$.
\begin{lemma}\label{lem-envelope}
Let $p>2$, let $(\Omega,\mu)$ be a probability space, and let
$X\subset L_p(\mu)$ be an $N$-dimensional real normed space. Suppose that
\[
\sup_{\|f\|\le1}\|f\|_{L_p(\mu)}^p\le M,
\]
and define
\[
s(x):=\sup_{\|f\|\le1}|f(x)|.
\]
Then $s$ is measurable and
\[
\int_\Omega s(x)^p\, \mu(dx)\le MN^{p/2}.
\]
\end{lemma}
\begin{proof}
Measurability follows by taking the supremum over a countable dense subset.
By Lewis's change of density \cite{Lewis1978} (see also \cite[Theorem~4.7]{KKLT}), there is a
probability density $\varrho$, positive on the essential support of $X$ and zero
elsewhere. With $T_\varrho f:=\varrho^{-1/p}f$ on $\{\varrho>0\}$, extended by zero to $\Omega$,
the space $T_\varrho(X)$ has an $L_2(\varrho\mu)$-orthonormal basis
$\phi_1,\ldots,\phi_N$ satisfying
\[
\sum_{j=1}^N|\phi_j(x)|^2=N
\quad\text{for }\varrho\mu\text{-almost every }x.
\]
Since $\varrho\mu$ is a probability measure and $p>2$, we have
\[
|T_\varrho f(x)|\le\sqrt N\|T_\varrho f\|_{L_2(\varrho\mu)}
\le\sqrt N\|T_\varrho f\|_{L_p(\varrho\mu)}
=\sqrt N\|f\|_{L_p(\mu)}.
\]
Thus $|f(x)|^p\le \varrho(x)N^{p/2}\|f\|_{L_p(\mu)}^p$
for $\mu$-almost every $x$. Taking the supremum over $\|f\|\le1$
and integrating gives the assertion.
\end{proof}

Considering the new norm $M^{1/p}\|\cdot\|$,
it suffices to prove the assertion
when $M=1$. Write $\xi:=(\xi^1,\ldots,\xi^m)$. For a random
variable $Y$ depending on $\xi$ and $r\ge1$, we will write
\[
\|Y\|_{L_r(\xi)}:=\bigl(\mathbb E_\xi|Y|^r\bigr)^{1/r}.
\]
Let
\begin{align*}
Z_\xi
&:=
\sup_{\|f\|\le 1}
\biggl|
\frac1m\sum_{j=1}^m|f(\xi^j)|^p
-\|f\|_{L_p(\mu)}^p
\biggr|
\quad\text{and}\quad
R_\xi
:=
\sup_{\|f\|\le1}
\biggl|
\frac1m\sum_{j=1}^m\eta_j|f(\xi^j)|^p
\biggr|,
\end{align*}
where $\eta_1,\ldots,\eta_m$ are independent Rademacher
random variables, independent of $\xi$.
Symmetrization gives
\begin{equation}\label{eq-symm}
\bigl(\mathbb E_\xi Z_\xi^r\bigr)^{1/r}
\le
2\bigl(\mathbb E_{\xi,\eta}R_\xi^r\bigr)^{1/r}
\quad\forall r\ge1.
\end{equation}
Let $s$ be defined as in Lemma \ref{lem-envelope}. For $k\ge0$, set
\[
A_k:=\{x\colon 2^{-(k+1)}H<s(x)^p\le2^{-k}H\},
\quad I_k:=\{j\colon \xi^j\in A_k\},\quad m_k:=|I_k|.
\]
Markov's inequality gives
\begin{equation}\label{eq-envelope-count}
\mathbb E_\xi m_k=m\mu(A_k)
\le\frac{2^{k+1}mN^{p/2}}H.
\end{equation}
For each nonempty $I_k$, equip the image
\[
E_k:=\{m^{-1/p}(f(\xi^j))_{j\in I_k}\colon f\in X\}
\]
with the quotient norm induced by $\|\cdot\|$ on $X$.
This norm is $p$-convex with the same constant $R$, and
\[
\sup_{x\in B_{E_k}}\|x\|_\infty^p\le\frac{2^{-k}H}m,
\qquad
\sup_{x\in B_{E_k}}\|x\|_p^p\le1+Z_\xi.
\]
Points where $s=0$ contribute nothing to $R_\xi$.
Applying Lemma~\ref{lem-signed-sum} on each nonempty level and using
Minkowski's inequality, we obtain, for every $r\ge1$,
\begin{align}\label{eq-conditional-estimate}
\bigl(\mathbb E_\eta R_\xi^r\bigr)^{1/r}
&\le C_1(p,R)\biggl[
\Bigl(\frac Hm\Bigr)^{1/p}(1+Z_\xi)^{1-1/p}
\\
&+\sqrt{\frac Hm}(1+Z_\xi)^{1/2}
\sum_{k=0}^\infty2^{-k/2}\bigl[\log_2(2(m_k+1))\bigr]^{1/2}\biggr] +C_2(p)\sqrt{\frac{rH}m}(1+Z_\xi)^{1/2},
\nonumber
\end{align}
Since $m_k$ is binomial, \cite[Corollary 1]{Ahle2022} gives
\[
\|m_k\|_{L_r(\xi)}\le\mathbb E_\xi m_k+\frac{r}{2}.
\]
Hence
\[
\bigl\|\log_2(2(m_k+1))\bigr\|_{L_r(\xi)}
\le\log_2\bigl(2(\mathbb E_\xi m_k+1)+r\bigr)+\frac1{\log2}.
\]
If $1\le r\le m/H$, then \eqref{eq-envelope-count} gives
\[
\bigl\|\log_2(2(m_k+1))\bigr\|_{L_r(\xi)}
\le C_3\log_2\Bigl(2+\frac{2^kmN^{p/2}}H\Bigr).
\]
Consequently, for this range of $r$,
\[
\sum_{k=0}^\infty2^{-k/2}\bigl\|\log_2(2(m_k+1))\bigr\|_{L_r(\xi)}^{1/2}
\le C_4\biggl[
\log_2\Bigl(2+\frac{mN^{p/2}}H\Bigr)
\biggr]^{1/2}.
\]
H\"older's inequality gives
\[
\bigl\|(1+Z_\xi)^{1/2}\bigl[\log_2(2(m_k+1))\bigr]^{1/2}\bigr\|_{L_r(\xi)}
\le(1+\|Z_\xi\|_{L_r(\xi)})^{1/2}
\bigl\|\log_2(2(m_k+1))\bigr\|_{L_r(\xi)}^{1/2},
\]
and, for $0<\theta<1$,
$\|(1+Z_\xi)^\theta\|_{L_r(\xi)}\le(1+\|Z_\xi\|_{L_r(\xi)})^\theta\le1+\|Z_\xi\|_{L_r(\xi)}$.
Taking the $L_r(\xi)$ norm of the conditional estimate \eqref{eq-conditional-estimate} and applying
symmetrization \eqref{eq-symm}, we therefore obtain, for $1\le r\le m/H$,
\begin{equation}\label{eq-sampling-moment}
\|Z_\xi\|_{L_r(\xi)}\le\biggl[
C_5(p,R)\biggl(
\Bigl(\frac Hm\Bigr)^{1/p}
+\sqrt{\frac Hm\log_2\Bigl(2+\frac{mN^{p/2}}H\Bigr)}
\biggr)+C_6(p)\sqrt{\frac{rH}m}
\biggr](1+\|Z_\xi\|_{L_r(\xi)}).
\end{equation}
Set
\[
m_0:=C(p,R)H\bigl(\varepsilon^{-p}
+\varepsilon^{-2}\log_2(2N)\bigr),
\]
where $C(p,R)\ge2$ will be chosen sufficiently large, and let $m\ge m_0$.
We have $\frac{H}{m}\le \frac{\varepsilon^p}{C(p,R)}$.
Since $m\mapsto m^{-1}\log_2(2+mN^{p/2}/H)$ is decreasing,
\begin{align*}
\frac Hm\log_2\Bigl(2+\frac{mN^{p/2}}H\Bigr)
&\le\frac H{m_0}\log_2\Bigl(2+\frac{m_0N^{p/2}}H\Bigr)\\
&\le\frac{C_7(p)\varepsilon^2}{C(p,R)}
\cdot\frac{\log_2 C(p,R)+p\log_2\varepsilon^{-1}+\log_2(2N)}
{\varepsilon^{-(p-2)}+\log_2(2N)}\le C_8(p)\frac{\log_2 C(p,R)}{C(p,R)}\varepsilon^2.
\end{align*}
Choose $c_0(p)>0$ sufficiently small and $C(p,R)$ sufficiently large,
and set
\[
r:=c_0(p)\frac{m\varepsilon^2}H.
\]
Then $1\le r\le m/H$ and, by \eqref{eq-sampling-moment},
\[
\|Z_\xi\|_{L_r(\xi)}\le \frac{\varepsilon}{4}(1+\|Z_\xi\|_{L_r(\xi)}).
\]
Since $\varepsilon<1$, it follows that $\|Z_\xi\|_{L_r(\xi)}\le\varepsilon/2$.
Markov's inequality now gives
\[
\mathbb P(Z_\xi>\varepsilon)
\le\Bigl(\frac{\|Z_\xi\|_{L_r(\xi)}}{\varepsilon}\Bigr)^r
\le2^{-r}
=\exp\Bigl(-c_0(p)(\log2)\frac{m\varepsilon^2}H\Bigr).
\]
Rescaling back to arbitrary $M$ and using homogeneity completes the proof.
\qed

\subsection{Proof of Corollary \ref{cor-sampling-Lp}}
Assume first that $X$ is real.
The $L_p(\mu)$-norm is $p$-convex with a constant depending
only on $p$. Applying Theorem~\ref{th-sampling-p-convex}
with this norm gives the assertion.

Now assume that $X$ is complex.
Let $\nu$ be the uniform probability measure on $[0,2\pi]$ and let
\[
\kappa_p:=
\int_0^{2\pi}|\cos t|^p\,\nu(dt)>0.
\]
We note that
\begin{equation}\label{eq-phase-identity}
\int_0^{2\pi}
\bigl|\operatorname{Re}(e^{\mathrm i(t+\theta)}w)\bigr|^p
\,\nu(dt)
=
\kappa_p|w|^p\quad
\forall w\in\mathbb C\quad \forall\theta\in\mathbb R.
\end{equation}

We consider the finite-dimensional real linear space
\[
Y:=\{g_f\colon f\in X\},
\qquad
g_f(\omega,t):=
\operatorname{Re}(e^{\mathrm it}f(\omega))
\]
of dimension $\dim Y= 2N$.
By \eqref{eq-phase-identity},
\[
\|g_f\|_{L_p(\mu\otimes\nu)}^p
=
\kappa_p\|f\|_{L_p(\mu)}^p.
\]
Moreover,
\[
\|g_f\|_\infty
\le
\|f\|_\infty
\le
\kappa_p^{-1/p}(KN)^{1/p}
\|g_f\|_{L_p(\mu\otimes\nu)}.
\]
Applying the real case to $Y$ with $K/\kappa_p$ in place of $K$ and enlarging $C(p)$ and decreasing $c(p)$ by factors
depending only on $p$,
for $m\ge C(p)KN\bigl(\varepsilon^{-p}
+\varepsilon^{-2}\log_2(2N)\bigr)$
and for independent samples
$(\xi^j,\theta_j)$ with distribution $\mu\otimes\nu$, with probability
at least
$1-\exp\bigl(-c(p)\frac{m\varepsilon^2}{KN}\bigr)$,
we have
\[
(1-\varepsilon)\|g_f\|_{L_p(\mu\otimes\nu)}^p
\le
\frac1m\sum_{j=1}^m
\bigl|\operatorname{Re}
(e^{\mathrm i\theta_j}f(\xi^j))\bigr|^p
\le
(1+\varepsilon)\|g_f\|_{L_p(\mu\otimes\nu)}^p
\quad\forall f\in X.
\]
In particular, for every $t\in [0, 2\pi]$,
\[
(1-\varepsilon)\|g_{e^{\mathrm it}f}\|_{L_p(\mu\otimes\nu)}^p
\le
\frac1m\sum_{j=1}^m
\bigl|\operatorname{Re}
(e^{\mathrm i(\theta_j+t)}f(\xi^j))\bigr|^p
\le
(1+\varepsilon)\|g_{e^{\mathrm it}f}\|_{L_p(\mu\otimes\nu)}^p
\quad\forall f\in X.
\]
Integrating the preceding inequalities with respect to $\nu$,
using \eqref{eq-phase-identity}, and cancelling $\kappa_p$ gives
\[
(1-\varepsilon)\|f\|_{L_p(\mu)}^p
\le
\frac1m\sum_{j=1}^m
\bigl|f(\xi^j)\bigr|^p
\le
(1+\varepsilon)\|f\|_{L_p(\mu)}^p
\quad\forall f\in X
\]
with probability
at least
$1-\exp\bigl(-c(p)\frac{m\varepsilon^2}{KN}\bigr)$.
This proves the complex case.
\qed

\section{Sums of entropy numbers of linear operators}
\label{sect-entropy}

For $1<q<2$ and for a sequence ${\bf Y}:=(Y_k)_{k=1}^\infty$ of Banach spaces, let
\[
\ell_q({\bf Y})
:=
\Bigl\{
({\bf y}_k)_{k=1}^\infty\colon
{\bf y}_k\in Y_k,\
\sum_{k=1}^\infty\|{\bf y}_k\|_{Y_k}^q<\infty
\Bigr\}.
\]
We note that $\ell_q({\bf Y})$, endowed with the norm
\[
\bigl\|({\bf y}_k)_{k=1}^\infty\bigr\|_{\ell_q({\bf Y})}:=\Bigl(\sum_{k=1}^\infty\|{\bf y_k}\|_{Y_k}^q\Bigr)^{1/q},
\]
is a Banach space. When all the spaces $Y_k$
coincide with one Banach space $Y$, we write
$\ell_q(Y)$.
For an element ${\bf v}\in \ell_q(Y)$
we define an operator
$L_{\bf v}\colon \ell_2\to Y$ by the formula
\[
L_{\bf v}{\bf x}:=\sum_{k=1}^\infty x_kv_k.
\]
Since
\[
\sum_{k=1}^\infty \|x_kv_k\|_Y
=
\sum_{k=1}^\infty |x_k|\|v_k\|_Y
\le
\Bigl(\sum_{k=1}^\infty \|v_k\|_Y^2\Bigr)^{1/2}
\|{\bf x}\|_{\ell_2}
\le
\|{\bf v}\|_{\ell_q(Y)}\|{\bf x}\|_{\ell_2},
\]
the operators $L_{\bf v}$ are well defined and bounded.

Theorem \ref{th-row-enrtopy} is a dual version of the following estimate.

\begin{theorem}\label{th-column-entropy}
Let $1<q<2$ and $T>0$.
There exists a constant $C:=C(q, T)>0$
such that, for every
Banach space $Y$ of type $q$
with constant $T$
and for every ${\bf v}\in \ell_q(Y)$,
one has
\[
\mathcal E_{2,q}(L_{\bf v})
\le C(q,T)\|{\bf v}\|_{\ell_q(Y)}.
\]
\end{theorem}

The following is a special case of
\cite[Theorem~4.1]{Cobos}, obtained by taking
$\varphi(t)=t^{1/q}$, $p=q$, and $v=2$.

\begin{theorem}[{\cite[Theorem~4.1]{Cobos}}]
\label{th-cobos}
Let $1<q<2$, let $Y$ be a Banach space of type $q$, and let
$J\colon\ell_1\to Y$ be a bounded linear operator.
For a sequence $\eta=(\eta_k)_{k=1}^\infty\in\ell_q$, let
\[
D_\eta\colon\ell_2\to\ell_1,
\quad
D_\eta {\bf x}=(\eta_kx_k)_{k=1}^\infty.
\]
Then
$\mathcal E_{2,q}(JD_\eta)<\infty$.
\end{theorem}

\begin{corollary}\label{cor-cobos}
Let $1<q<2$ and let $Y$ be a Banach space of type $q$. Then
\[
\mathcal E_{2,q}(L_{\bf v})<\infty\quad \forall {\bf v}\in \ell_q(Y).
\]
\end{corollary}

\begin{proof}
Let $\eta_k: = \|v_k\|_Y$ and let $\eta:=(\eta_k)_{k=1}^\infty$. Then $\eta \in \ell_q$ and
\[
L_{\bf v}{\bf x} = JD_\eta {\bf x},
\]
where $J\colon \ell_1\to Y$,
\[
J{\bf x}:= \sum_{\substack{k\ge1\\v_k\ne0}} x_k \frac{v_k}{\|v_k\|_Y}.
\]
The assertion now follows from Theorem \ref{th-cobos}.
\end{proof}

\begin{lemma}\label{lem-column}
Let $1<q<2$ and let $Y$ be a Banach space of type $q$. Then
\[
C(Y):=\sup_{\|{\bf v}\|_{\ell_q(Y)}=1}\mathcal E_{2,q}(L_{\bf v})
<\infty.
\]
\end{lemma}

\begin{proof}
Assume that the assertion is not true.
Then there is a sequence ${\bf v}_n=(v_{n,k})_{k=1}^\infty\in \ell_q(Y)$
such that
\[
\|{\bf v}_n\|_{\ell_q(Y)}=1\quad\textnormal{and}\quad
\mathcal E_{2,q}(L_{{\bf v}_n})\ge2^{2n}.
\]
Let $r\mapsto (n(r), k(r))$ be any bijection from $\mathbb{N}$ to $\mathbb{N}^2$
and let ${\bf v}:= (2^{-n(r)}v_{n(r), k(r)})_{r=1}^\infty$.
Then ${\bf v}\in \ell_q(Y)$ and
\[
\varepsilon_j(L_{\bf v})\ge 2^{-n}
\varepsilon_j(L_{{\bf v}_n}),
\]
implying that
\[
\mathcal E_{2,q}(L_{\bf v})\ge 2^n\quad \forall n\in \mathbb{N},
\]
which contradicts Corollary \ref{cor-cobos}.
\end{proof}

\begin{proof}[Proof of Theorem~\ref{th-column-entropy}]
Let
\[
C(q,T):=
\sup\{C(Y)\colon Y\text{ is of type }q
\text{ with constant }T\}.
\]
Suppose that $C(q,T)=+\infty$.
Then there is a sequence of Banach spaces $Y_n$ of type $q$
with constant $T$ such that $C(Y_n)\to+\infty$.
For the collection ${\bf Y}:=(Y_n)_{n=1}^\infty$,
it is readily verified that $\ell_q({\bf Y})$
is also of type $q$ with constant $T$.

Fix $n\in\mathbb N$.
For $(v_k)_{k=1}^\infty\in\ell_q(Y_n)$, let
${\bf v}_k\in\ell_q({\bf Y})$ have $n$th coordinate $v_k$
and all other coordinates zero, and let
${\bf v}:=({\bf v}_k)_{k=1}^\infty\in \ell_q(\ell_q({\bf Y}))$.
Then
\[
\|{\bf v}\|_{\ell_q(\ell_q({\bf Y}))}
=
\|(v_k)_{k=1}^\infty\|_{\ell_q(Y_n)}
\]
and
\[
\varepsilon_j(L_{\bf v})
\ge
\varepsilon_j(L_{(v_k)_{k=1}^\infty})
\quad\forall j\in\mathbb N.
\]
Hence,
\[
\mathcal E_{2,q}(L_{\bf v})\ge \mathcal E_{2,q}(L_{(v_k)_{k=1}^\infty})
\]
and
therefore
\[
C(\ell_q({\bf Y}))\ge C(Y_n)
\quad\forall n\in\mathbb N.
\]
Since $C(Y_n)\to+\infty$, this contradicts
Lemma~\ref{lem-column}.
\end{proof}

\begin{proof}[Proof of Theorem \ref{th-row-enrtopy}]
By the duality properties of the entropy numbers, see, for instance, the proof of   \cite[Theorem~4.3.7]{AGM2015} or \cite[Corollary~2(i) and Remarks~(1)--(2)]{BPST1989}, we have
\[
\mathcal{E}_{2, q}(A)\le C(q)
\mathcal{E}_{2, q}(A^*).
\]
Now,
$A^*\colon \ell_2^m\to X^*$,
$A^*z = \sum_{k=1}^m z_ka_k$.
By Theorem \ref{th-column-entropy},
\[
\mathcal{E}_{2, q}(A^*)\le
C(q,T)
\biggl(\sum_{k=1}^m\|a_k\|_{X^*}^q\biggr)^{1/q},
\]
which completes the proof.
\end{proof}

\section{Sharpness examples}
\label{sect-examples}

Throughout this section, $p>2$ is fixed.

\begin{example}\label{ex-1}
Let $K\ge1$, $N\ge 4$, and consider
\[
\Omega=\{0,1,\ldots,N\},
\qquad
\mu(\{j\})=\frac{1}{KN}\quad \forall j\in \{1, \ldots, N\},
\qquad
\mu(\{0\})=1-\frac {1}{K}.
\]
The space
\[
X:=\operatorname{span}
\bigl\{\mathbf1_{\{j\}}\colon 1\le j\le N\bigr\}
\]
has dimension $N$ and satisfies
\[
\|f\|_\infty\le (KN)^{1/p}\|f\|_{L_p(\mu)}
\quad\forall f\in X.
\]

Let $\xi^1,\ldots,\xi^m$ be independent samples with
distribution $\mu$, and let $M$ be the number of atoms from
$\{1,\ldots,N\}$ absent from the sample.
It is clear that
\[
\mathbb{P}\biggl((1-\varepsilon)\|f\|_{L_p(\mu)}^p
\le
\frac1m\sum_{j=1}^m|f(\xi^j)|^p
\le
(1+\varepsilon)\|f\|_{L_p(\mu)}^p
\quad\forall f\in X
\biggr)\le \mathbb{P}(M=0).
\]
Since
\[
M=\sum_{k=1}^N{\bf 1}_{A_k}, \quad A_k:=\{k\not\in\{\xi^1, \ldots, \xi^m\}\},
\]
we have
\[
\mathbb EM=N\Bigl(1-\frac{1}{KN}\Bigr)^m.
\]
Since, for $k\ne k'$,
\[
\mathbb{E}[{\bf 1}_{A_k}{\bf 1}_{A_{k'}}] =
\mathbb{E}[{\bf 1}_{A_k\cap A_{k'}}] =
\Bigl(1-\frac{2}{KN}\Bigr)^m
\le
\Bigl(1-\frac{1}{KN}\Bigr)^{2m} = \mathbb{E}[{\bf 1}_{A_k}]\mathbb{E}[{\bf 1}_{A_{k'}}],
\]
we have
$\operatorname{Var}(M)\le\mathbb EM$, and
Chebyshev's inequality gives
\[
\mathbb P(M=0) \le \mathbb P\bigl((M-\mathbb{E}M)^2\ge (\mathbb{E}M)^2\bigr)
\le
\frac{\operatorname{Var}(M)}{(\mathbb EM)^2}
\le\frac1{N(1-\frac{1}{KN})^m}.
\]
Now, if
$m\le\frac {KN}4\log N$,
then
\[
\Bigl(1-\frac1{KN}\Bigr)^m
\ge e^{-2m/(KN)}\ge N^{-1/2}.
\]
Thus,
\[
\mathbb{P}\biggl((1-\varepsilon)\|f\|_{L_p(\mu)}^p
\le
\frac1m\sum_{j=1}^m|f(\xi^j)|^p
\le
(1+\varepsilon)\|f\|_{L_p(\mu)}^p
\quad\forall f\in X
\biggr)\le N^{-1/2}.
\]
In particular, for the discretization inequality to hold with
probability greater than $1/2$, it is necessary that
$m>\frac {KN}4\log N$.
\end{example}

\begin{example}\label{ex-2}
Let $\Omega=\{1,\ldots,N\}^2$ and let $\mu$ be the uniform probability
measure on $\Omega$. Let
\[
X:=\{f\colon\Omega\to\mathbb R\colon
f(k,j)=a_k+b_j,\ a,b\in\mathbb R^N\}.
\]
Since the linear map $(a,b)\mapsto(a_k+b_j)_{k,j=1}^N$
from $\mathbb R^{2N}$ to $\mathbb R^{N^2}$ has image $X$
and a one-dimensional kernel, we have $\dim X=2N-1$. Averaging the identity
\[
f(k,j)=f(k,t)+f(s,j)-f(s,t)
\]
over $s,t$ gives
\[
\|f\|_\infty^p
\le3^{p}(2N+1)\|f\|_{L_p(\mu)}^p
\le3^{p+1}N\|f\|_{L_p(\mu)}^p.
\]
Thus, we may take $K=3^{p+1}$.

Suppose that nonnegative weights $\lambda_{kj}$ satisfy
\begin{equation}\label{eq-example-discr}
(1-\varepsilon)\|f\|_{L_p(\mu)}^p\le
\sum_{k,j=1}^N\lambda_{kj}|f(k,j)|^p
\le(1+\varepsilon)\|f\|_{L_p(\mu)}^p
\quad\forall f\in X,
\end{equation}
where $0<\varepsilon\le1/4$.
We estimate the number of positive weights.

Fix $k_0\in\{1, \ldots, N\}$, and let
\[
r_{k_0}:=\sum_{j=1}^N\lambda_{k_0j}
\quad\text{and}\quad
S_{k_0}:=\{j\colon \lambda_{k_0j}>0\}.
\]
Testing the inequality \eqref{eq-example-discr}
on the function ${\bf 1}_{\{k=k_0\}}\in X$ gives
\begin{equation}\label{eq-example-r-est}
r_{k_0}\ge\frac{1-\varepsilon}{N},
\end{equation}
so $|S_{k_0}|\ge1$.
If $|S_{k_0}|\le N/2$, let $t:=|S_{k_0}|^{-1/p}$ and consider
\[
f_\pm(k,j):=
\mathbf1_{\{k=k_0\}}\pm t\mathbf1_{\{j\in S_{k_0}\}}\in X.
\]
For
\[
\Delta:=(1+t)^p-(1-t)^p\ge2pt,
\]
we obtain
\[
\sum_{k,j=1}^N\lambda_{kj}|f_+(k,j)|^p
-
\sum_{k,j=1}^N\lambda_{kj}|f_-(k,j)|^p
=r_{k_0}\Delta\quad\text{and}\quad
\|f_+\|_{L_p(\mu)}^p-\|f_-\|_{L_p(\mu)}^p
=\frac{|S_{k_0}|\Delta}{N^2}.
\]
Applying the upper bound in \eqref{eq-example-discr} to $f_+$
and the lower bound to $f_-$ gives
\[
\Bigl(r_{k_0}-\frac{|S_{k_0}|}{N^2}\Bigr)\Delta
\le
\varepsilon
\bigl(\|f_+\|_{L_p}^p+\|f_-\|_{L_p}^p\bigr).
\]
Since $|S_{k_0}|\le N/2$ and $\varepsilon\in (0, 1/4]$, from \eqref{eq-example-r-est}, we obtain
\[
r_{k_0}-\frac{|S_{k_0}|}{N^2}\ge\frac1{4N}.
\]
Moreover, since $|S_{k_0}|t^p=1$ and $0<t\le1$,
\[
\|f_+\|_{L_p(\mu)}^p+\|f_-\|_{L_p(\mu)}^p
=
\frac{
|S_{k_0}|\bigl((1+t)^p+(1-t)^p\bigr)
+2(N-|S_{k_0}|)+2(N-1)|S_{k_0}|t^p}{N^2}
\le\frac{2^{p-1}+4}{N}.
\]
Consequently,
\[
\frac{p}{2N}|S_{k_0}|^{-1/p}
\le\frac{\varepsilon(2^{p-1}+4)}N,
\]
and hence
\[
|S_{k_0}|\ge c_1(p)\varepsilon^{-p},
\quad
c_1(p):=
\Bigl(\frac{p}{2(2^{p-1}+4)}\Bigr)^p.
\]
Thus,
\[
\#\{(k,j)\colon\lambda_{kj}>0\}
\ge N\min\{N/2,c_1(p)\varepsilon^{-p}\}
\ge c_2(p) K\dim X\min\{\dim X,\varepsilon^{-p}\}.
\]
Repeated sampling points can be combined into a single weight,
so the same lower bound applies to the number of sampling points.

In particular, for sufficiently large $N$, taking
$\varepsilon=(\dim X)^{-1/p}$ shows that sampling discretization
of the $L_p(\mu)$ norm on $X$ requires at least $c(p)(\dim X)^2$ sampling
points. For any fixed $q<p$ and $\alpha\ge0$,
\[
\dim X(\log(2\dim X))^\alpha\varepsilon^{-q}
=
(\dim X)^{1+q/p}(\log(2\dim X))^\alpha
=o((\dim X)^2).
\]
Therefore, a sample-size bound of this order cannot guarantee
sampling discretization in general.
\end{example}

\begin{example}
Let $K\ge2$, and let $A$ be a measurable set with
$\mu(A)=1/K$. For the one-dimensional space
\[
X:=\operatorname{span}\{\mathbf1_A\},
\]
we have
\[
\|f\|_\infty\le K^{1/p}\|f\|_{L_p(\mu)}
\quad\forall f\in X.
\]
For independent samples $\xi^1,\ldots,\xi^m$
with the distribution $\mu$, let
\[
M:=\sum_{j=1}^m\mathbf1_A(\xi^j),
\quad
a:=\frac1K,
\quad
\lambda:=ma.
\]
Then $M$ has the binomial distribution with parameters $m,a$,
and
\[
\mathbb{P}\biggl((1-\varepsilon)\|f\|_{L_p(\mu)}^p
\le
\frac1m\sum_{j=1}^m|f(\xi^j)|^p
\le
(1+\varepsilon)\|f\|_{L_p(\mu)}^p
\quad\forall f\in X
\biggr)\le \mathbb P\bigl(|M-\lambda|\le \varepsilon\lambda\bigr).
\]
Assume that $0<\varepsilon\le1/4$ and $m\ge \varepsilon^{-2}K$.
We show that
\[
\mathbb P\bigl(|M-\lambda|>\varepsilon\lambda\bigr)
\ge C\exp\Bigl(-c\frac{\varepsilon^2 m}{K}\Bigr).
\]

For $u\in(0,1)$, let
\[
h_a(u):=
u\log\frac ua+(1-u)\log\frac{1-u}{1-a}.
\]
Since $\log v\le v-1$, we have
\[
h_a(u)
\le
u\Bigl(\frac ua-1\Bigr)
+(1-u)\Bigl(\frac{1-u}{1-a}-1\Bigr)
=
\frac{(u-a)^2}{a(1-a)}.
\]
We apply \cite[Theorem~9]{ZhangZhou2020} with
$\beta=2$ and $x=2\varepsilon\lambda$.
Its assumptions hold because
\[
0<x\le\frac{\lambda}{2}
\le\frac{m(1-a)}2
\quad\text{and}\quad
x+\lambda\ge1.
\]
Consequently, for absolute constants $c_1,C>0$,
\[
\mathbb P(M-\lambda\ge2\varepsilon\lambda)
\ge
C\exp\bigl(-c_1m h_a((1+2\varepsilon)a)\bigr).
\]
Moreover,
\[
m h_a((1+2\varepsilon)a)
\le
\frac{4ma\varepsilon^2}{1-a}
\le 8\frac{m\varepsilon^2}{K}.
\]
Therefore, for another absolute constant $c_2>0$,
\[
\mathbb P\bigl(|M-\lambda|>\varepsilon\lambda\bigr)
\ge
\mathbb P(M-\lambda\ge2\varepsilon\lambda)
\ge
C\exp\Bigl(-c_2\frac{m\varepsilon^2}{K}\Bigr).
\]
Hence the exponential rate $m\varepsilon^2/K$ in the
failure probability in Corollary \ref{cor-sampling-Lp} is optimal up to absolute constant factors.
\end{example}

\section*{Use of AI Tools}

This paper is a follow-up to our initial arXiv
preprint \cite{DKM} and was prepared after we consulted ChatGPT
in early October 2026 about the problems studied there.
The use of Cobos's result \cite{Cobos}, suggested by ChatGPT,
allowed us to remove the extra factor $(\log\log N)^{p-1}$
from our sampling discretization result presented in \cite{DKM}.
Using results from \cite{Kosov2021}, we simplified the initial
argument suggested by ChatGPT.

The authors acknowledge the use of ChatGPT during the preparation
of the manuscript to assist with the development of the main
mathematical argument and with language and stylistic editing.
All AI-generated material used was reviewed and revised by the authors.
The authors take full responsibility for the content of the paper.

\end{document}